\documentclass[12pt]{article}
\usepackage{pgf,tikz}
\usepackage{calc,graphicx}
\usepackage{graphics}
\everymath{\displaystyle}
\usepackage{graphicx}
\usepackage{color}
\usepackage{amssymb}
\usepackage{amsmath}
\newtheorem{prethm}{{\bf Theorem}}

\newenvironment{thm}{\begin{prethm}{\hspace{-0.5
				em}{\bf .}}}{\end{prethm}}
\newtheorem{prelemma}{{\bf Lemma}}

\newtheorem{preex}{{\bf Example}}

\newtheorem{preprop}{{\bf Proposition}}

\newenvironment{prop}{\begin{preprop}{\hspace{-0.5em}{\bf .}}}{\end{preprop}}
\newtheorem{precor}{{\bf Corollary}}

\newenvironment{cor}{\begin{precor}{\hspace{-0.5
				em}{\bf .}}}{\end{precor}}
\newtheorem{preremark}{{\bf Remark}}

\newtheorem{preprob}{{\bf Problem}}

\newtheorem{predefin}{{\bf Definition}}

\newtheorem{preconj}{{\bf Conjecture}}

\newtheorem{preprobb}{{\bf Problem}}

\newtheorem{prelem}{{\bf Theorem}}

\newtheorem{precla}{{\bf Claim}}

\newenvironment{cla}{\begin{precla}{\hspace{-0.5
em}{\bf .}}}{\end{precla}}

\newenvironment{proof}{{\bf Proof.}\rm }{\hfill{$\Box$}}

\newtheorem{presolution}{{\bf Solution.}}

\def\newpic#1{}
\def\qed{\ifhmode\unskip\nobreak\fi\quad\ifmmode\Box\else$\Box$\fi}

\title{\vspace{-2.2cm}\Large\bf\noindent More results on the $z$-chromatic number of graphs}
\author{\large\bf Abbas Khaleghi~~~~~~Manouchehr Zaker\footnote{mzaker@iasbs.ac.ir}
\vspace{5mm}\\
Department of Mathematics,\\
Institute for Advanced Studies in Basic Sciences,\\
Zanjan 45137-66731, Iran\\
}

\date{\bf Submitted to Discrete Appl. Math. on September 08, 2022}

\begin{document}
\maketitle

\begin{abstract}
\noindent By a $z$-coloring of a graph $G$ we mean any proper vertex coloring consisting of the color classes $C_1,
\ldots, C_k$ such that $(i)$ for any two colors $i$ and $j$ with $1 \leq i < j \leq k$, any vertex of color $j$ is adjacent to a vertex of color $i$, $(ii)$ there exists a set $\{u_1, \ldots, u_k\}$ of vertices of $G$ such that $u_j \in C_j$ for any $j \in \{1, \ldots, k\}$ and $u_k$ is adjacent to $u_j$ for each $1 \leq j \leq k$ with $j \not=k$, and $(iii)$ for each $i$ and $j$ with $i \not= j$, the vertex $u_j$ has a neighbor in $C_i$. Denote by $z(G)$ the maximum number of colors used in any $z$-coloring of $G$. Denote the Grundy and {\rm b}-chromatic number of $G$ by $\Gamma(G)$ and ${\rm b}(G)$, respectively.  The $z$-coloring is an improvement over both the Grundy and b-coloring of graphs. We prove that $z(G)$ is much better than $\min\{\Gamma(G), {\rm b}(G)\}$ for infinitely many graphs $G$ by obtaining an infinite sequence $\{G_n\}_{n=3}^{\infty}$ of graphs such that $z(G_n)=n$ but $\Gamma(G_n)={\rm b}(G_n)=2n-1$ for each $n\geq 3$. We show that acyclic graphs are $z$-monotonic and $z$-continuous. Then it is proved that to decide whether $z(G)=\Delta(G)+1$ is $NP$-complete even for bipartite graphs $G$. We finally prove that to recognize graphs $G$ satisfying $z(G)=\chi(G)$ is $coNP$-complete, improving a previous result for the Grundy number.
\end{abstract}

\noindent {\bf Keywords:} Graph coloring; First-Fit coloring; Grundy number; {\rm b}-chromatic number; $z$-chromatic number; $z$-coloring

\noindent {\bf AMS Classification:} 05C15, 05C85

\section{Introduction}

\noindent All graphs in this paper are undirected without any loops and multiple edges. In a graph $G$, $\Delta(G)$ denotes the maximum degree of $G$. Let $v$ be a vertex in $G$, $N(v)$ and $N[v]$ denote the set of neighbors of $v$ and the closed neighborhood of $v$ in $G$, respectively. Also for any subset $S$ of vertices in $G$, by $G[S]$ we mean the subgraph of $G$ induced by the elements of $S$. A complete graph on $n$ vertices is denoted by $K_n$. The union of two vertex disjoint graphs $G_1$ and $G_2$ is the graph $G_1 \cup G_2$ with vertex set $V(G_1)\cup V(G_2)$ and edge set $E(G_1)\cup E(G_2)$. The join of two disjoint graphs $G_1$ and $G_2$ is the graph $G_1 \oplus G_2$ obtained from $G_1 \cup G_2$ by putting an edge between any vertex in $G_1$ and any vertex in $G_2$. A proper vertex coloring of a graph $G$ is an assignment of colors $1, 2, \ldots$ to the vertices of $G$ such that any two adjacent vertices receive distinct colors. By a color class we mean a subset of vertices having a same color. The smallest number of colors used in a proper coloring of $G$ is called the chromatic number of $G$ and is denoted by $\chi(G)$. A proper edge coloring of $G$ is defined similarly. The  minimum number of distinct colors required for a proper edge coloring of $G$ is denoted by $\chi'(G)$. It was proved in \cite{3d} that it is $NP$-complete to decide whether a given 3-regular graph $G$ is $3$-edge colorable. We refer to \cite{BM} for the terminology not defined here. The Grundy and b-coloring are two well-known techniques for proper coloring of graphs to be defined as follow.

\noindent By a Grundy-coloring of a graph $G$ we mean any proper vertex coloring of $G$ consisting of color classes say $C_1, \ldots, C_k$ such that for each $i < j$ any vertex in $C_j$ has a neighbor in $C_i$. The Grundy number (also called the First-Fit chromatic number) of a graph $G$, denoted by $\Gamma(G)$ (also by $\chi_{\sf FF}(G)$) is the maximum number of colors used in any Grundy-coloring of $G$. Clearly, $\Gamma(G)\leq \Delta(G)+1$. The literature is full of papers concerning the Grundy number and First-Fit coloring of graphs e.g. \cite{GL, HS, Z1, Z2}. The $NP$-completeness of determining the Grundy number was proved for the complement of bipartite graphs in \cite{Z1} and \cite{Z2} and for bipartite graphs in \cite{HS}. Graphs satisfying $\Gamma(G)=\chi(G)$ are called well-colored graphs in \cite{Z2}, where it was proved that the recognition of well-colored graphs is $coNP$-complete.

\noindent A {\rm b}-coloring (or color-dominating coloring) of a graph $G$ is a proper vertex coloring in which any color class contains a vertex (color-dominating vertex) adjacent to at least one vertex in every other color class. The {\rm b}-chromatic number {\rm b}$(G)$ is the largest integer $k$ such that there is a {\rm b}-coloring of $G$ using $k$ colors. Clearly, ${\rm b}(G)\leq \Delta(G)+1$. To determine ${\rm b}(G)$ is $NP$-complete but has a polynomial-time solution for trees \cite{Ir}. A graph $G$ is called b-continuous in \cite{Fai} if for any integer $k$, $\chi(G) \leq k \leq {\rm b}(G)$, $G$ admits a b-coloring with $k$ colors. In \cite{Bono} Bonomo et al. introduced the concept of {\rm b}-monotonicity. A graph $G$ is {\rm b}-monotonic if ${\rm b}(H_1)\leq {\rm b}(H_2)$ for every induced subgraph $H_1$ of $G$ and every induced subgraph $H_2$ of $H_1$.

\noindent We note that the Grundy and b-coloring are two color-suppressing techniques by applying certain recoloring operations. Extensions of the recoloring operation used in the b-coloring was studied in \cite{PR}. A sophisticated recoloring operation and corresponding color-suppressing technique is introduced in \cite{Z5} under the name of $z$-coloring. A $z$-coloring of a graph $G$ is any proper vertex coloring consisting of the color classes $C_1,
\ldots, C_k$ such that $(i)$ for any two colors $i$ and $j$ with $1 \leq i < j \leq k$, any vertex of color $j$ is adjacent to a vertex of color $i$, $(ii)$ there exists a set $\{u_1, \ldots, u_k\}$ of vertices of $G$ such that $u_j \in C_j$ for any $j \in \{1, \ldots, k\}$ and $u_k$ is adjacent to $u_j$ for each $1 \leq j \leq k$ with $j \not=k$, and $(iii)$ for each $i$ and $j$ with $i \not= j$, the vertex $u_j$ has a neighbor in $C_i$. It was proved in \cite{Z5} that any graph admits a $z$-coloring by an efficient procedure. Denote by $z(G)$ the maximum number of colors used in any $z$-coloring of $G$. In a $z$-coloring $C$ of a graph $G$ using $k$ colors, a vertex $v$ is called nice vertex if $v$ has color $k$ in $C$ and $v$ is adjacent to at least $k-1$ color-dominating vertices with $k-1$ distinct colors in $C$. Note that $z$-coloring is obtained from a recoloring technique until we get a Grundy-coloring with a nice vertex in the graph. It is then an improvement over the Grundy and color-dominating colorings. We have $z(G)\leq \min\{\Gamma(G),{\rm b}(G)\} \leq \Delta(G)+1$.  We say that a graph $G$ is $z$-continuous if and only if for any $k$, $\chi(G)\leq k \leq z(G)$, there exists a $z$-coloring of $G$ using $k$ colors. Also a graph $G$ is $z$-monotonic if $z(H_{1})\geq z(H_{2})$ for every induced subgraph $H_{1}$ of $G$ and every induced subgraph $H_{2}$ of $H_{1}$. These two properties do not hold in general. The graph $K_{n,n}\setminus nK_2$ is not $z$-continuous and $K_{n,n}\setminus (n-1)K_2$ is not $z$-monotonic.

\noindent {\bf The outline of the paper is as follows.} In Proposition \ref{ju} we obtain a relation for the $z$-chromatic number of the join and union of graphs. Then in Proposition \ref{example} we present a sequence of graphs $\{G_n\}_{n=3}^{\infty}$ such that for each $n$, $z(G_n)< \Gamma(G_n)={\rm b}(G_n)$. We show in the rest of Section 2 that acyclic graphs are $z$-continuous and $z$-monotonic. In Section 3 we prove that to decide whether $z(G)=\Delta(G)+1$ is $NP$-complete for bipartite graphs $G$. It is finally proved that to recognize graphs $G$ satisfying $z(G)=\chi(G)$ is $coNP$-complete.

\section{Some general results}

\noindent In this section we first obtain relations for the $z$-chromatic number of $G_1\oplus G_2$ and $G_1\cup G_2$ in terms of $z(G_1)$ and $z(G_2)$. These relations will be used in a later result. Similar relations for the Grundy and b-chromatic number were obtained in \cite{Bli}. The relations are similar except that if $G_1$ and $G_2$ are two vertex-disjoint graphs then ${\rm b}(G_1\cup G_2)\geq max \{{\rm b}(G_1),{\rm b}(G_2)\}$ and the inequality is strict in many cases. But for the $z$-chromatic number we have equality.

\begin{prop}\label{ju}
Let $G_1$ and $G_2$ be two vertex-disjoint graphs. Then

\noindent (i) $z(G_1\oplus G_2)=z(G_1) + z(G_2)$.

\noindent (ii) $z(G_1\cup G_2)=max \{z(G_1), z(G_2)\}$.
\end{prop}

\noindent \begin{proof}
\noindent (i) Set $z(G_1)=k_1$, $z(G_2)=k_2$, $G=G_1\oplus G_2$ and assume that $k_1\geq k_2$. Let $c_1$ be a $z$-coloring of $G_1$ with $k_1$ colors and $c_2$ be a $z$-coloring of $G_2$ with $k_2$ colors. Define a $z$-coloring $c$ of $G$ with $k_1+k_2$ colors by letting
\begin{equation*}
c(v)=
\begin{cases}
c_1(v)+k_2 & v\in V(G_1)
\\
c_2(v) & v\in V(G_2)
\end{cases}
\end{equation*}
\noindent It is easy to check that if $v\in G_1$ is a nice vertex in $c_1$ then $v$ is a nice vertex in $z$-coloring $c$ of $G$ with $k_1+k_2$ colors. Hence, $z(G)\geq z(G_1)+z(G_2)$.

\noindent Now, assume on the contrary that $z(G)> k_1+k_2$. Then for some integer $t\geq 1$, $z(G)=k_1+k_2+t$. Consider a $z$-coloring $c$ of $G$ with $k_1+k_2+t$ colors. Since no color can appear in both $G_1$ and $G_2$, then either for $i=1$ or for $i=2$, there are strictly more than $k_i$ colors in $G_i$. Assume that it happens for $i=1$. Let $C_1, C_2, \ldots, C_j$ be these color classes in $G_1$ with $j>k_1$. Define a $z$-coloring $c'$ of $G_1$ as follows. For every vertex $v\in G_1$ if $c(v)\in C_i$, $1\leq i\leq j$, then set $c'(v)=i$. The proper coloring $c'$ is a $z$-coloring of $G_1$ with $j$ colors since $c$ is a $z$-coloring of $G$. This contradicts with the maximality of $z(G_1)$.

\noindent (ii) Set $z(G_1\cup G_2)=k$. Clearly, $k\geq max \{z(G_1), z(G_2)\}$. Assume that $k> max \{z(G_1), z(G_2)\}$. Consider a $z$-coloring $c$ of $G_1\cup G_2$ with $k$ colors. Let $w$ be a nice vertex in $c$ of color $k$. Then there exist $k-1$ color-dominating vertices of distinct colors which are adjacent to $w$. Since $G_1$ and $G_2$ are vertex disjoint, we may assume that $w$ and all of these color-dominating vertices belongs to say $G_1$. Let $c_1$ be the coloring of $G_1$ by restricting $c$ to $G_1$. It follows that $c_1$ is a $z$-coloring of $G_1$. But $c_1$ has more colors than $z(G_1)$, a contradiction.
\end{proof}

\noindent The next result shows that for infinitely many graphs $z(G)$ is much better than $\min \{\Gamma(G),{\rm b}(G)\}$. An infinite sequence of graphs $\{H_n\}_{n=1}^{\infty}$ was constructed in \cite{Z5} such that $\min \{\Gamma(H_n),{\rm b}(H_n)\} \rightarrow \infty$ as $n\rightarrow \infty$ but $z(H_n)\leq 3$, for each $n$. Proposition \ref{example} provides another sequence $G_n$ such that $\min \{\Gamma(G_n),{\rm b}(G_n)\} - z(G_n) \rightarrow \infty$.

\begin{prop}
There exists an infinite sequence of graphs $\{G_n\}_{n=3}^{\infty}$ such that for each $n$, $\Gamma(G_n)={\rm b}(G_n)=2n-1$ but $z(G_n)=\chi(G_n)=n$.\label{example}
\end{prop}

\noindent \begin{proof}
For any positive integer $n\geq 3$, we construct a graph $G_n$ as follows. First, consider a complete graph $K_n$ on a vertex set $\{v_1,\ldots,v_n\}$. Then to each vertex $v_i$ attach a complete graph $K(v_i)$ on $n$ vertices $\{v_{1}^i,\ldots, v_{n-1}^i,v_n^i\}$ such that $v_i$ and $v_n^i$ from $K(v_i)$ are identified, i.e. $v_i=v_n^i$. We have $V(K(v_i))\cap V(K(v_j))=\emptyset$, for every $i$, $j$ with $1\leq i < j \leq n$. Finally, for each $i$, $1\leq i\leq n$, attach $n-1$ leaves to $v_1^i$ of $K(v_i)$. Denote the resulting graph by $G_n$. Note that $\Delta(G_n)=2n-2$ and $\chi(G_n)=n$. A general form of the graph $G_n$ is depicted in Figure \ref{zg2}.

\begin{figure}[h]
\centering
\includegraphics[scale=0.50]{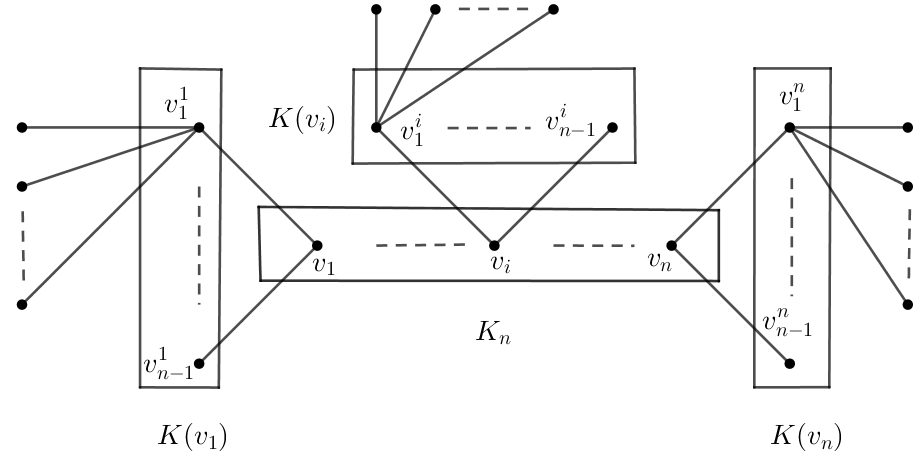}
\caption{The graph $G_n$.}
\label{zg2}
\end{figure}

\noindent We show that $\Gamma(G_n)=2n-1$. By $\Delta(G_n)=2n-2$ and the known upper bound $\Gamma(G)\leq \Delta(G)+1$ we have $\Gamma(G_n)\leq 2n-1$. Now we present a partial Grundy-coloring for $G_n$ using $2n-1$ colors. Assign colors $n, n+1, \ldots, 2n-1$ to the vertices $v_1,\ldots,v_n$, respectively. Then for any vertex other than $v_n^i$ in $K(v_i)$, $i\in \{1,\ldots,n\}$, assign a color $c$ from the set $\{1, \ldots, n-1\}$, such that no pair of vertices in $K(v_i)$ receive a same color. We can easily extend this partial Grundy-coloring to a Grundy-coloring of $G_n$ using $2n-1$ colors. Figure \ref{zg} illustrates an small instance of this situation.

\begin{figure}[h]
\centering
\includegraphics[scale=0.50]{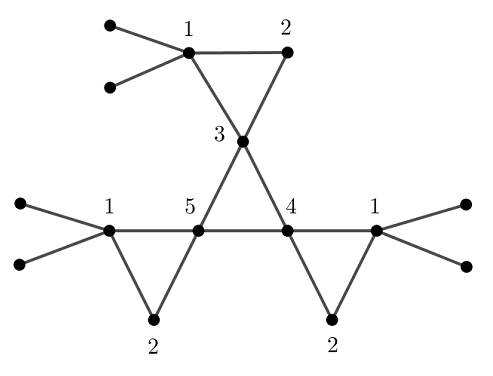}
\caption{A partial Grundy-coloring of $G_3$ with $5$ colors.}
\label{zg}
\end{figure}

\noindent We prove that ${\rm b}(G_n)=2n-1$. We present a partial vertex coloring for $D=\{v_1, \ldots, v_n, v_1^1, v_1^2, \ldots, v_1^{n-1}\}$ using $2n-1$ distinct colors. Assign colors $n, \ldots, 2n-1$ to $v_1, \ldots, v_n$, respectively. Then assign colors $1, \ldots, n-1$ to $v_1^1, \ldots, v_1^{n-1}$, respectively. Figure \ref{zg1} illustrates this partial coloring for $G_4$ in which the elements of $D$ are displayed as gray vertices. Now, using the leaf vertices we can extend this partial coloring to a {\rm b}-coloring of $G_n$ using $2n-1$ colors, where the vertices of $D$ are the color-dominating vertices. Hence, ${\rm b}(G_n)\geq 2n-1$. Now, the inequality ${\rm b}(G_n)\leq \Delta(G_n)+1=2n-1$ proves the required equality.

\begin{figure}[h]
\centering
\includegraphics[scale=0.50]{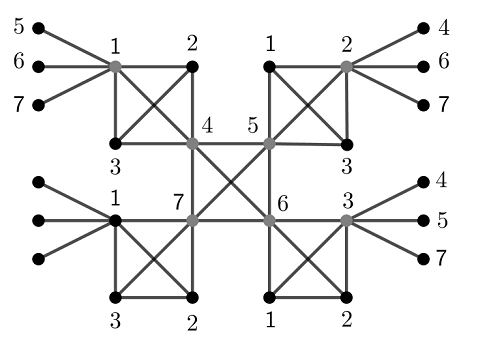}
\caption{A partial coloring of $G_4$ with $7$ colors and color-dominating gray vertices}
\label{zg1}
\end{figure}

\noindent To complete the proof we argue that $z(G_n)= n$ for each $n\geq 3$. Assume on the contrary that $G_n$ admits a $z$-coloring using $n+1$ or more colors. Let $v$ be a nice vertex of color $n+1$ in a $z$-coloring of $G_n$. The vertex
$v$ needs at least $n$ neighbors of degree at least $n$. By the construction of $G_n$, $v$ should be in $\{v_1, \ldots, v_n\}$. Without loss of generality, assume that $v=v_1$. Since the only neighbors of $v_1$ having degree at least $n$ are
the $n$ vertices $v_2, \ldots, v_n$ and $v_1^1$ from $K(v_1)$, then we may assume that $v_1^1$ is a color-dominating vertex of color $1$. Hence, each color in $\{2, \ldots, n\}$ should appear in $N(v_1^1)$. It follows that there exists $w\in K(v_1)$ with a color $j \in \{2, \ldots, n\}$. Then $w$ does not have any neighbors of color $i$ with $i<j$, a contradiction with the first property of $z$-coloring. Hence $z(G_n)\leq n$. Now, $\chi(G_n)\leq z(G_n)$ implies $z(G_n)=n$.
\end{proof}

\noindent In the following we prove that acyclic graphs are $z$-continuous and $z$-monotonic.

\begin{prop}
Any acyclic graph is $z$-continuous.
\end{prop}

\noindent \begin{proof}
By Proposition \ref{ju} $(ii)$, it is enough to prove the proposition for trees. Let $T$ be a tree which admits a $z$-coloring using $k$ colors, where $k \geq 3$. We obtain a $z$-coloring for $T$ using exactly $k-1$ colors. As explained before, $T$ contains the tree $R_k$ as induced subtree. By the construction of $R_k$, $R_{k-1}$ is contained in $R_k$ and hence in $T$. Consider a $z$-coloring for $R_{k-1}$ using $k-1$ colors. This partial $z$-coloring of $T$ is easily extended to a $z$-coloring of whole $T$ with exactly $k-1$ colors.
\end{proof}

\begin{prop}
Any acyclic graph is $z$-monotonic.\label{monotone}
\end{prop}

\noindent \begin{proof}
By Proposition \ref{ju} $(ii)$, it is enough to prove for trees $T$. We show that for any vertex $v$ of $T$, $z(T-v)\leq z(T)$. Set $z(T-v)=k$. Let $T_{1},\ldots,T_{m}$ be components of $T-v$, we have $z(T-v)=max \{z(T_{1}),\ldots,z(T_{m})\}$, without loss of generality, assume that $z(T-v)=z(T_{1})=k$, so we have a $z$-coloring for $T_{1}$ with $k$ colors. We extend it to a $z$-coloring for $T$, since the vertex $v$ is adjacent to only one vertex of $T_{1}$, it is enough to color the vertex $v$ and each of $T_{i}$'s, $(i\neq 1)$ with colors $1$ or $2$. We obtain a $z$-coloring for $T$ with $k$ colors.
\end{proof}

\noindent We need the concept of $k$-atoms. For each positive integer $k$, a class of graphs denoted by ${\mathcal{A}}_k$ was constructed in \cite{Z2} which satisfies the following property. The Grundy number of any graph $G$ is at least $k$ if and only if $G$ contains an induced subgraph isomorphic to some element of ${\mathcal{A}}_k$. Any graph in ${\mathcal{A}}_k$ is called {\it $k$-atom}. The concept of atom graphs has been also used in study of b-coloring of graphs \cite{EGT}. For any positive integer $k$, there exists exactly one tree $k$-atom, denoted by $T_{k}$. For $k=1,2$, $T_k$ is isomorphic to the complete graph on one and two vertices, respectively. Assume that $T_k$ is constructed for $k\geq 2$, then $T_{k+1}$ is obtained from $T_k$ by attaching one leaf to each vertex of $T_k$ so that $|V(T_{k+1})|=2|V(T_k)|$. Proposition \ref{ztk} determines $z(T_k)$. It was proved in \cite{Z5} that there exists a unique tree $R_k$ such that $z(R_k)=k$ and for every tree $T$, $z(T)\geq k$ if and only if $T$ contains a subtree isomorphic to $R_k$. 

\begin{prop}\label{ztk}
For any integer $k$, $z(T_k)=\lceil (k+1)/2 \rceil$.
\end{prop}

\noindent \begin{proof}
By the construction of $T_k$, its degree sequence is:
$$k-1, k-1, k-2, k-2, k-3, k-3, k-3, k-3, k-4, \ldots, 2, \ldots, 2, 1, \ldots,1$$
\noindent For any $i\in\{1, \ldots, k-2\}$, there are exactly $2^{k-1-i}$ vertices of degree $i$ in the degree sequence. To prove $z(T_k)\geq \lceil (k+1)/2 \rceil$, it's enough by Proposition \ref{monotone} to obtain a subgraph $H$ of $T_k$ with $z(H)\geq \lceil (k+1)/2 \rceil$. In a Grundy-coloring $c$ of $T_k$ with $k$ colors, there exists one vertex say $v_k$ of color $k$. Considering $v_k$ as a root, it has $k-1$ children $v_1, v_2, \ldots, v_{k-1}$ of colors $1, 2, \ldots, k-1$, respectively. Namely, $c(v_i)=i$ for each $i$. Add $v_k, v_{k-1}, v_{\lceil k/2 \rceil}$ to $H$ and assign new colors $\lceil (k+1)/2 \rceil, \lceil (k+1)/2 \rceil-1, \ldots, 2, 1$ to $v_k, v_{k-1}, \ldots, v_{\lceil k/2 \rceil +1}, v_{\lceil k/2 \rceil}$, respectively. Note that $\lceil (k+1)/2 \rceil = k- \lceil k/2 \rceil +1$. we extend this partial coloring and also the subgraph $H$ itself so that $v_i$ is a color-dominating vertex of color $i-\lceil k/2 \rceil +1$, for each $i\in \{\lceil k/2 \rceil, \ldots, k\}$. This makes $v_k$ to be a nice vertex in $H$.
Consider the vertex $v_{\lceil k/2 \rceil}$ whose color in $c$ is $\lceil k/2 \rceil$ and then has $\lceil k/2 \rceil-1$ children of colors $1, 2, \ldots, \lceil k/2 \rceil-1$. We add the children of colors $2, \ldots, \lceil k/2 \rceil-1$ and their descendant vertices with their colors from $c$ to $H$. Then $v_{\lceil k/2 \rceil}$ becomes color-dominating of color $1$ in $H$. The technique of extension corresponding to other vertices in $v_{\lceil k/2 \rceil +1}, \ldots, v_k$ is similar. We obtain $z(T_k)\geq z(H) \geq \lceil (k+1)/2 \rceil$.

\noindent To prove $z(T_k)\leq \lceil (k+1)/2 \rceil$, assume on the contrary that $T_k$ admits a $z$-coloring $c'$ using $\lceil (k+1)/2 \rceil+1$ colors. Let $v$ be a nice vertex of color $\lceil (k+1)/2 \rceil+1$ in $c'$. Vertex $v$ needs at least $\lceil (k+1)/2 \rceil$ neighbors of degree at least $\lceil (k+1)/2 \rceil$. But it can be easily proved by an induction on $k$ that any vertex of degree at least $t$ in $T_k$ has at most $k-t$ neighbors of degree at least $t$ in $T_k$. In particular, any vertex of degree at least $\lceil (k+1)/2 \rceil$ has at most $k-\lceil (k+1)/2 \rceil = \lfloor (k-1)/2 \rfloor$ neighbors of degree at least $\lceil (k+1)/2 \rceil$. Hence such a vertex $v$ does not exist.
\end{proof}

\noindent The following result is a corollary of Propositions \ref{ztk} and \ref{monotone}.

\begin{cor}
For any tree $T$, $\Gamma(T)\leq 2z(T)-1$ and equality holds for the tree atoms.
\end{cor}

\noindent \begin{proof}
Set $\Gamma(T)=k$. As explained earlier, $T$ contains $T_k$ as subgraph, then $z(T_k)=\lceil (k+1)/2 \rceil \leq z(T)$.
\end{proof}

\section{Complexity results}

\begin{figure}[h]
\centering
\includegraphics[scale=0.60]{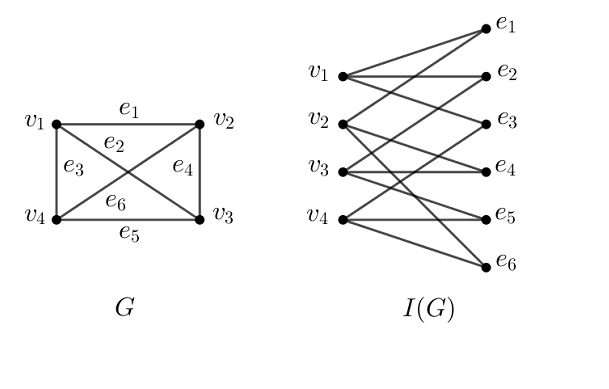}
\caption{A 3-regular graph $G$ and its incidence graph $I(G)$.}
\label{K4}
\end{figure}

\noindent In this section we show that determining the $z$-chromatic number is $NP$-complete even for bipartite graphs.
In the following by $M_{n,n}$ we mean the graph obtained from the complete bipartite graph $K_{n,n}$ by removing a perfect matching, i.e. $M_{n,n}=K_{n,n}\setminus nK_2$. For any graph $G$, the vertex-edge incidence graph of $G$, denoted by $I(G)$ is the bipartite graph with the bipartition $V(I(G))=V(G) \cup E(G)$ in which an arbitrary edge $e=uv$ of $G$ (as a vertex in $I(G)$) is adjacent to its two endvertices $u, v$ in $I(G)$. An example of $I(G)$ is depicted in Figure \ref{K4}.

\begin{thm}\label{bipar}
It is $NP$-complete to decide if a bipartite graph $H$ satisfies $z(H)=\Delta(H)+1$.\label{thm2}
\end{thm}

\noindent \begin{proof}
The problem belongs to $NP$ because let $C$ be a proper vertex coloring of $H$. We can determine in ${\mathcal{O}}(|V(H)|+|E(H)|)$ steps whether it is a $z$-coloring with at least $k$ colors. To show that the problem is $NP$-hard, we present a reduction from $3$-edge-colorability of $3$-regular graphs, which is $NP$-complete as mentioned before.

\noindent Assume that $A$ and $B$ are the bipartite sets in $M_{n+2,n+2}=K_{n+2,n+2}\setminus(n+2)K_{2}$. Let $F$ be the graph presented in Figure \ref{F} in which a vertex is distinguished as the vertex $f$. We construct from $F$ and $M_{n+2,n+2}$ a new graph $T$ as follows, we connect vertex $f$ of $F$ to $n-1$ vertices of $B$. For any positive integer $i$, let $T_{i}$ be an isomorphic copy of $T$ in which the vertices of $M_{n+2,n+2}$ are labelled $\{v_{1}^i,\ldots,v_{i-1}^{i},v_{i+1}^{i},\ldots, v_{n}^i, x^i, y^i, q^i\}$ and $\{u_{1}^i,\ldots,u_{i-1}^{i},u_{i+1}^{i},\ldots u_{n}^i, a^i, b^i, c^i\}$. Figure \ref{T1} illustrates $T_{i}$.

\begin{figure}[h]
\centering
\includegraphics[scale=0.60]{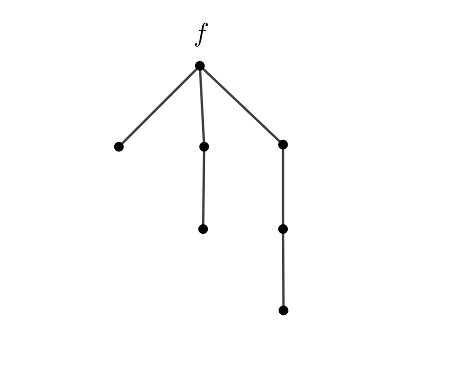}
\caption{The graph $F$.}
\label{F}
\end{figure}

\begin{figure}[h]
\centering
\includegraphics[scale=0.65]{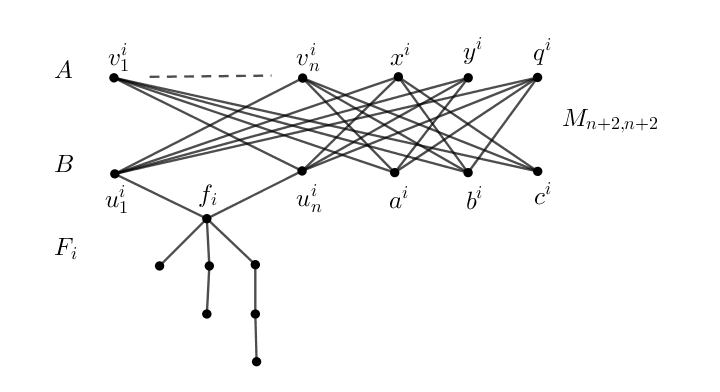}
\caption{The graph $T_i$.}
\label{T1}
\end{figure}

\noindent Let $G$ be any $3$-regular graph with $n$ vertices. Set $V(G)=\{v_1,\ldots, v_n\}$ and $E(G)=\{e_{1},\ldots,e_{m}\}$. Let $I(G)$ be the vertex-edge incidence graph of $G$. We construct gradually from $I(G)$ a new graph $H$ as follows. At first let $H$ be the $I(G)$ itself. Then for each vertex $e_{i} \in E(G)$ (as a vertex of $I(G)$), add to $H$ a copy $M_{3,3}(e_{i})$ isomorphic to $M_{3,3}$ and identify one vertex of $M_{3,3}(e_i)$ with the vertex $e_{i}$ of $I(G)$. In other words $e_{i}\in M_{3,3}(e_{i})$. Recall that $M_{n,n}=K_{n,n}\setminus nK_{2}$. Next, add a new vertex $w$ to $H$ adjacent to all the vertices of $V(G)$ in $I(G)$ and add copies $M_{n+3,n+3}^{1}$, $M_{n+3,n+3}^{2}$ and $M_{n+3,n+3}^{3}$ all isomorphic to $M_{n+3,n+3}$ to $H$. Then choose arbitrary vertices $v_{n+1}$, $v_{n+2}$ and $v_{n+3}$, respectively from each copy and add the edges $wv_{n+1}$, $wv_{n+2}$ and $wv_{n+3}$. Finally, for any $i\in \{1, \ldots, n\}$, add the graph $T_i$ to $H$ and put an edge between $v_i$ and all vertices in $\{v_1^i, \ldots, v_{i-1}^{i}, v_{i+1}^{i}, \ldots, v_{n}^i\}\subseteq V(T_i)$. An illustration of $H$ is presented in Figure \ref{HH}. We have
$$V(H)=V(G) \cup \bigcup_{i=1}^m V(M_{3,3}(e_i)) \cup (\bigcup_{i=1}^n V(T_i)) \cup (\bigcup_{i=1}^3 V(M_{n+3,n+3}^i)) \cup\{w\}.$$

\begin{figure}[h]
\centering
\includegraphics[scale=0.55]{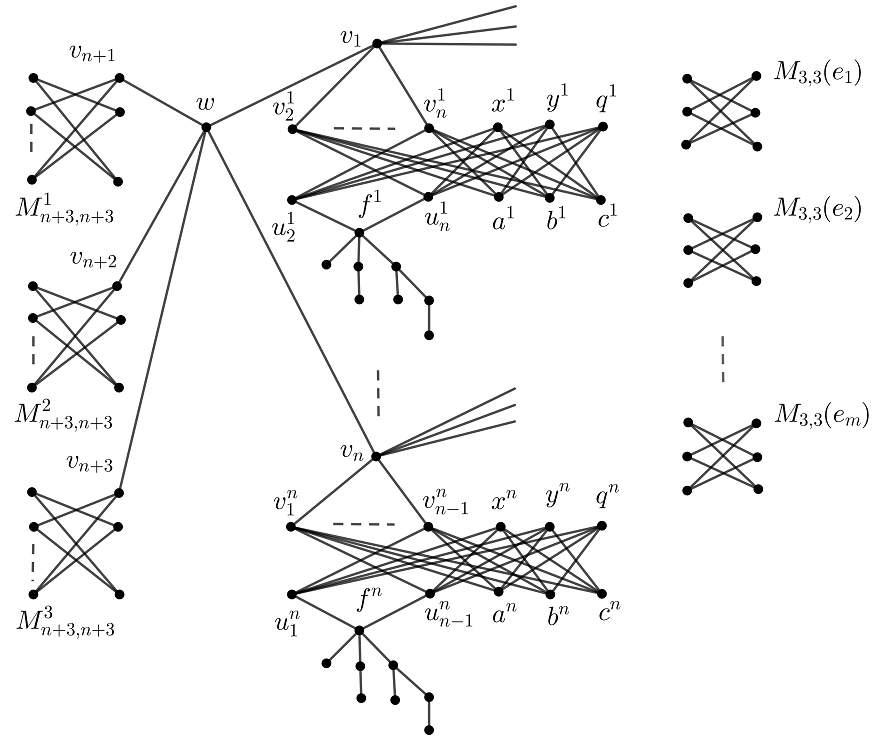}
\caption{The graph $H$.}
\label{HH}
\end{figure}

\noindent We have the following facts concerning $H$.

\noindent $(i)$ $d_{H}(w)=n+3$, and $d_{H}(v_{i})=n+3$, for $1 \leq i \leq n$.

\noindent $(ii)$ $d_{H}(v_{j}^{i})=n+2$, for $1\leq i \leq n$ and $1\leq j \leq n$, $j\neq i$, since a vertex in $M_{n+2,n+2}$ of $T_{i}$ has degree $n+1$ and also $v_{i}$ is adjacent to $v_{j}^{i}$.

\noindent $(iii)$  $d_{H}(u_{j}^{i})=n+2$, for $1\leq i \leq n$ and $1\leq j \leq n$, $j\neq i$, since a vertex in $M_{n+2,n+2}$ of $T_{i}$ has degree $n+1$ and also $f_{i}$ is adjacent to $u_{j}^{i}$.

\noindent $(iv)$ $d_{H}(e_{j})=4$, for $1\leq j \leq m$, since $e_{j}$ has two neighbors in $I$ and two in $M_{3,3}(e_{i})$.

\noindent $(v)$ $ d_{H}(x^{i})=d_{H}(y^{i})=d_{H}(q^{i})=d_{H}(a^{i})=d_{H}(b^{i})=d_{H}(c^{i})=n+1$ and $d(f^i)=n+2$, for $1 \leq i \leq n$ and $d_{H}(v_{n+1})=d_{H}(v_{n+2})=d_{H}(v_{n+3})=n+3$.

\noindent $(vi)$ $\Delta(H)=n+3$ and the only vertices with degree $n+3$ are $w$ and $v_{i}$, for $1\leq i \leq n+3$.

\noindent To prove the theorem we show that $z(H)=\Delta(H)+1=n+4$ if and only if $G$ is $3$-edge-colorable.

\noindent Assume first that $c$ is a $z$-coloring using $n+4$ colors for $H$. We prove that the edges of $G$ can be properly colored using $3$ colors.

\noindent Any nice vertex of color $n+4$ in $c$ needs at least $n+3$ neighbors of degree at least $n+3$. The only vertex having these properties is $w$, therefore $c(w)=n+4$. Since $d_{H}(w)=n+3$, $w$ has exactly one neighbor colored $i$, for every $1 \leq i \leq n+3$. Then for every $1\leq i \leq n+3$, $c(v_i)\in \{1,2,\ldots,n+3\}$. Since $\{v_1, \ldots, v_n, v_{n+1}, v_{n+2}, v_{n+3}\}$ are color-dominating, then for any vertex $v$ from the latter set and for any $j \in \{1, 2, \ldots, n+3\}$ with $j\neq c(v)$, $v$ should have a neighbor of color $j$.

\begin{cla}\label{cla1}
For any $i \in \{1, \ldots, n\}$, $c(v_i)\in \{4,5,\ldots,n+3\}$.
\end{cla}

\noindent {\bf Proof of Claim 1:}

\noindent Assume on contrary that there exists an $i\in\{1,2,\ldots,n\}$ such that $c(v_i)\in \{1,2,3\}$.
Without loss of generality, suppose that $c(v_1)=1$. Since $v_1$ is color-dominating of degree $n+3$ and color 1, then each color in $\{2,3,\ldots,n+3\}$ should appear in $N[v_1]$. Vertex $v_1$ has exactly $n-1$ neighbors in $T_1$. The other neighbors of $v_1$ are of type $e_p$ for some $p\in \{1, \ldots, m\}$. Note that by the property $(iv)$, $d_{H}(e_p)=4$ for every $p$ and then $c(e_p)\leq 5$. It follows that for some $l,r,s \in \{1,\ldots, m\}$, $v_1$ should be adjacent to $e_l$, $e_s$ and $e_r$, respectively from $M_{3,3}(e_l)$, $M_{3,3}(e_s)$ and $M_{3,3}(e_r)$ such that $c(e_l),c(e_s),c(e_r)\in \{2,3,4,5\}$. Now, one of the following cases holds.

\noindent {\bf Case 1.} $c(e_{l}),c(e_{s}),c(e_{r}) \in \{2,3,4\}$

\noindent In this case the neighbors of $v_{1}$ with colors $\{5,6,\ldots,n+3\}$ should appear in the neighborhood of $v_1$ in the part $A$ from $M_{n+2,n+2}$ of $T_1$. It implies that the vertices of $A$ which are adjacent to $v_{1}$ must have one neighbor of each color $2,3,4$. These neighbors have to be in the $B$ part of $T_1$. Therefore $c(a^1),c(b^1),c(c^1)\in \{2,3,4\}$.

\noindent {\bf Case 2.} $c(e_{l}),c(e_{s}),c(e_{r}) \in \{2,3,5\}$

\noindent In this case the neighbors of $v_{1}$ with colors $\{4,6,\ldots,n+3\}$ should appear in the neighborhood of
$v_1$ in the part $A$ from $M_{n+2,n+2}$ of $T_1$. It implies that the vertices of $A$ which are adjacent to $v_{1}$ must have one neighbor of each color $2,3,5$. These neighbors have to be in the $B$ part of $T_1$. Therefore $c(a^1),c(b^1),c(c^1)\in \{2,3,5\}$.

\noindent {\bf Case 3.} $c(e_{l}),c(e_{s}),c(e_{r}) \in \{2,4,5\}$

\noindent A similar argument proves that in this case $c(a^{1}),c(b^{1}),c(c^{1})\in \{2,4,5\}$.

\noindent {\bf Case 4.} $c(e_{l}),c(e_{s}),c(e_{r}) \in \{3,4,5\}$

\noindent A similar argument proves that in this case $c(a^{1}),c(b^{1}),c(c^{1})\in \{3,4,5\}$. These four cases are illustrated in Figures \ref{T1T2} and \ref{T3T4}.

\begin{figure}[h]
\centering
\includegraphics[scale=0.50]{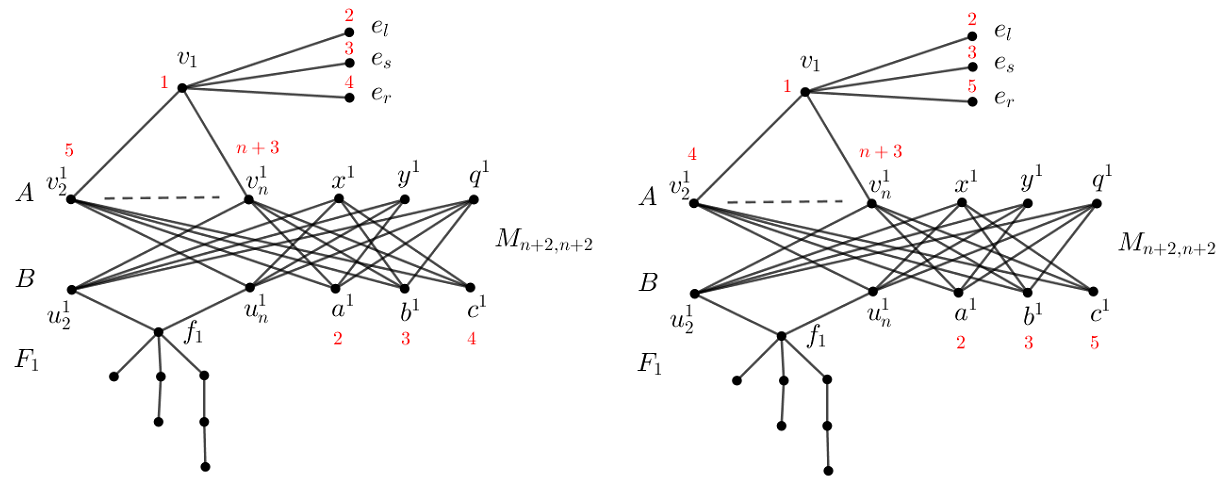}
\caption{A coloring of $T_{1}$ such that $c(e_{l}),c(e_{s}),c(e_{r}) \in \{2,3,4\}$ and $\{2,3,5\}$.}
\label{T1T2}
\end{figure}

\begin{figure}[h]
\centering
\includegraphics[scale=0.50]{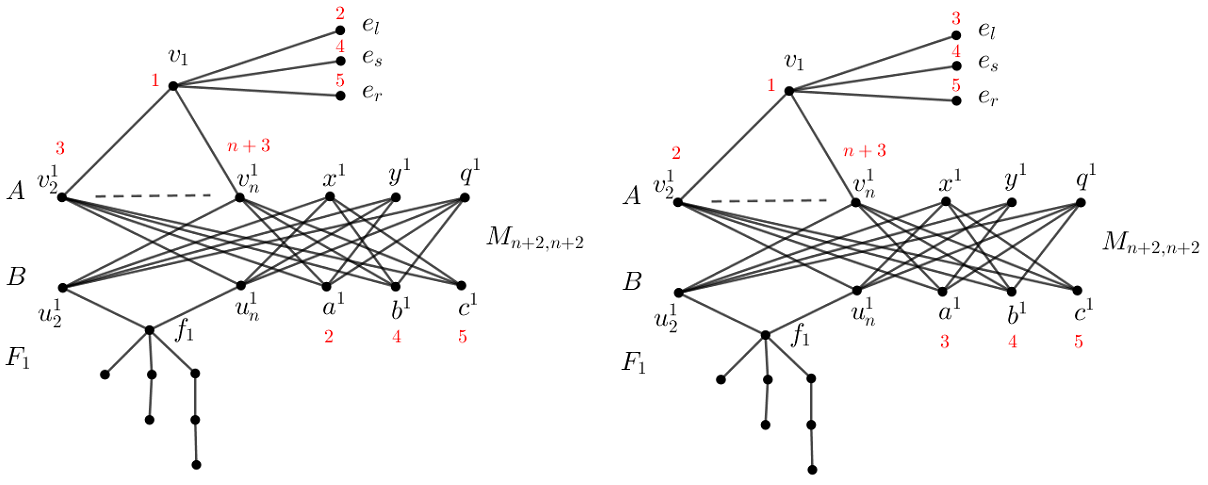}
\caption{A coloring of $T_{1}$ such that $c(e_{l}),c(e_{s}),c(e_{r}) \in \{2,4,5\}$ and $\{3,4,5\}$.}
\label{T3T4}
\end{figure}

\noindent Note that in the all cases $c(u_{j}^{i})=c(v_{j}^{i})$ for $1\leq i \leq n$ and $1\leq j \leq n-1$.
\noindent In each of the four cases, the vertices $a^1$, $b^1$ and $c^1$ don't have any neighbor with color 1, that is a contradiction. With the same argument, it follows that for every $1\leq i \leq n$, $c(v_i)\notin \{2,3\}$.
This completes the proof of Claim 1.

\begin{cla}\label{cla2}
For every $1\leq j\leq m$, $c(e_j)\in \{1,2,3\}$.
\end{cla}

\noindent {\bf Proof of Claim 2:}

\noindent Let $v_i$ be an arbitrary vertex with $1\leq i \leq n$. Vertex $v_i$ is color-dominating of degree $n+3$ and by Claim \ref{cla1}, $c(v_i)\in \{4,\ldots,n+3\}$. Hence, for each $t\in \{1, 2, 3\}$, $v_i$ has exactly one neighbor of color $t$. If for one neighbor say $e_j\in M_{3,3}(e_j)$ of $v_i$ we have $c(e_j)\notin \{1,2,3\}$, then by the first property of $z$-coloring, there is some $p\in \{1,\ldots,n\}$ with $p \neq i$ such that $v_p$ is adjacent to $e_j$ and $c(v_p)=3$. Note that the vertex $v_p$ is not appeared in $M_{3,3}(e_j)$, that is a contradiction by claim \ref{cla1}. The situation is depicted in Figure \ref{v_p}.

\begin{figure}[h]
\centering
\includegraphics[scale=0.50]{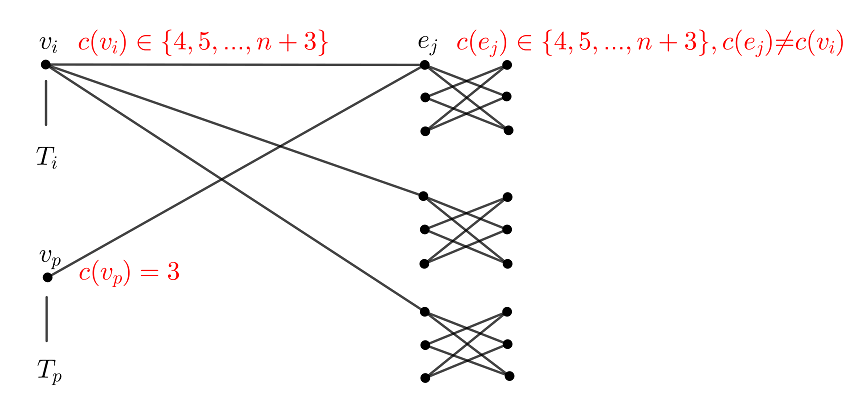}
\caption{$v_p$, $T_i$ and $T_{p}$ as in the proof of Claim \ref{cla2}.}
\label{v_p}
\end{figure}

\noindent Finally, $w$ should have three color-dominating neighbors of each color from $\{1,2,3\}$. Since $N_H(w)=\{v_1, \ldots, v_n, v_{n+1}, v_{n+2}, v_{n+3}\}$, then $c(v_{n+1}), c(v_{n+2}), c(v_{n+3}) \in \{1,2,3\}$.

\noindent We now prove that $c$ induces a proper 3-edge-coloring $c'$ of $G$.

\noindent As proved before $c(w)=n+4$. Without loss of generality, we may assume that for every $1\leq i \leq n$, $c(v_i)=i+3$. Since $v_i$ is adjacent to $w$ which has a color greater than $n+3$ and $d(v_i)=n+3$, there are only $n+2$ vertices remaining for the other $n+2$ colors. Therefore for some $l,s,r \in\{1,\ldots,m\}$ there exist three vertices adjacent to $v_i$ as $e_l$, $e_s$ and $e_r$ respectively from $M_{3,3}(e_l)$, $M_{3,3}(e_s)$ and $M_{3,3}(e_r)$ with different colors, by Claim \ref{cla2}, these colors can only be 1, 2, 3. Furthermore, for any $j\in\{1,\ldots,m\}$, $e_j$ has at most two neighbors in $H$ having colors at most $3$. Define $c'(e_j)=c(e_j) \in \{1, 2, 3\}$. Indeed, the three edges incident to $v_i$ in $G$, which are also adjacent to $v_i$ in $H$ takes different colors.

\noindent We show the properness of $c'$. If $c'(e_j)=c'(e_t)$ such that $e_j$ and $e_t$ has a common endvertex $v_i$, then $v_i$ can not be color-dominating vertex. Hence, $c'$ is a proper $3$-edge coloring.

\noindent Assume now that $G$ has a $3$-edge-coloring $\theta$ using the colors $1,2,3$. We obtain a $z$-coloring $c$ of $H$ with $n+4$ colors as follow. First, for any $1 \leq j \leq m$, greedily color the vertices of $M_{3,3}(e_j)$ such that $e_j$ is colored $\theta(e_j)$ and define $c(e_j)=\theta (e_j)$. By doing this, every vertex in $V(G)\subseteq V(H)$ has one neighbor of each color in $\{1, 2, 3\}$.

\noindent Next, for every $1 \leq i \leq n$, greedily color $v_i$ and $T_i$ such that $c(v_i)=c(f^i)=i+3$ and $c(a^i)=c(x^i)=1$, $c(b^i)=c(y^i)=2$, $c(c^i)=c(q^i)=3$. Also for any $1 \leq j \leq n$ with $j\neq i$, set $c(v_{j}^{i})=c(u_{j}^{i})=j+3$.

\noindent Next, greedily color $M_{n+3,n+3}^1$, $M_{n+3,n+3}^2$ and $M_{n+3,n+3}^3$ with $\{1, \ldots, n+3\}$ in such a way that $w$ has one neighbor colored $i$, for any $i \in \{1, 2, 3\}$. This means that $v_{n+1}$, $v_{n+2}$ and $v_{n+3}$ are color-dominating vertices with colors 1, 2 and 3, respectively.

\noindent Finally, assign $n+4$ to $w$. Clearly for any $i\in \{1,\ldots, n+3\}$, $w$ has a color-dominating neighbor of each color $i$. It follows that $z(H)=n+4$, since $\Delta(H)=n+3$.
\end{proof}

\noindent The following construction and corresponding result shall be used in the next result.

\begin{prop}\label{prop1}
\noindent Let $P_4$ be the path on four vertices in which the two vertices of degree $1$ are $v_1$ and $v_4$. Let $G$ be a connected bipartite graph of maximum degree $\Delta(G)$. Construct a graph $H$ from $G$, $P_4$ and $K_{\Delta(G)}$ by joining $v_1$ to all vertices of $G$ and joining $v_4$ to all vertices of $K_{\Delta(G)}$. Then
$$z(H)=max \{z(G), \Delta(G)\}+1$$.
\end{prop}

\begin{figure}[h]
\centering
\includegraphics[scale=0.55]{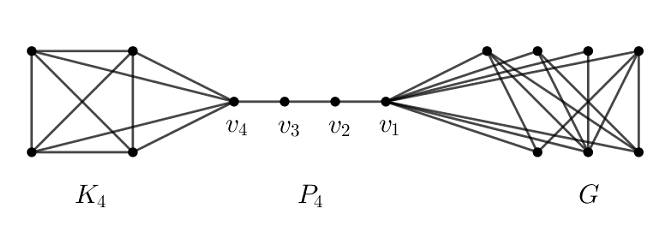}
\caption{The construction of $H$ in Proposition \ref{prop1} for a bipartite graph $G$ with $\Delta(G)=4$.}
\label{figureH}
\end{figure}

\noindent \begin{proof}
The graph $H$ corresponding to a bipartite graph $G$ with $\Delta(G)=4$ is illustrated in Figure \ref{figureH}. It is clear that $z(H)\geq max \{z(G), \Delta(G)\}+1$. We prove that $z(H)\leq max \{z(G), \Delta(G)\}+1$. First suppose that $z(G)\geq \Delta(G)$. Assume on the contrary that $z(H) > \max \{z(G), \Delta(G)\} + 1$. Then for some $t \geq 1$, we have $z(H)=z(G)+t+1$. So there is a $z$-coloring with $z(G)+t+1$ colors for $H$. Let $v$ be a nice vertex of color $z(G)+t+1$ in this $z$-coloring. The vertex $v$ needs at least $z(G)+t$ neighbors of degree at least $z(G)+t$. Since $d_{H}(v)=\Delta(G)$ for each $v\in K_{\Delta(G)} \oplus v_4, v\neq v_4$, then $v\notin V(K_{\Delta(G)}\oplus v_{4})$. Also $v\not\in \{v_2, v_3\}$ because the degree of $v_2$ and $v_3$ is two. Therefore $v\in G\oplus v_{1}$. If $v=v_1$ then none of $v_1$'s color-dominating neighbors can be $v_2$ because the degree of $v_2$ is two. Since the vertex $v_1$ is joined to $G$ and $z(G\oplus v_{1})=z(G)+1$. It follows that every $z$-coloring of $H$ using more than $z(G)+1$ colors reduces to a $z$-coloring with more than $z(G)+1$ colors in $G\oplus v_{1}$, that is a contradiction. A similar argument proves the case of $z(G)\leq \Delta(G)$.
\end{proof}

\noindent It was proved in \cite{Z2} that to recognize graphs $G$ satisfying $\Gamma(G)=\chi(G)$ is $coNP$-complete. We generalize this result to the $z$-chromatic number.

\begin{thm}
The following decision problem is $coNP$-complete.

\noindent {Instance}: graph $G$.

\noindent {Question}: $z(G)=\chi(G)?$\label{prob1}
\end{thm}

\noindent \begin{proof}
First observe that $z(G)>\chi(G)$ if and only if there exists two $z$-colorings for $G$ using say $k$ and $k'$ colors such that $k>k'$. Hence, to prove that the problem belongs to $coNP$, a short certificate is a pair $(C,C')$ of $z$-colorings for $G$ such that $C$ uses strictly more colors than $C'$. To check that an assignment of colors to the vertices of $G$ is a $z$-coloring is done in ${\mathcal{O}}(|E(G)|)$ time steps. It remains to prove that the problem is $NP$-hard. By Theorem \ref{bipar} it is $NP$-complete to decide whether $z(G) =\Delta (G)+1$ for a given bipartite graph $G$. We introduce a polynomial time reduction from the latter problem to the complement our problem. We transform a given bipartite graph $G$ of maximum degree $\Delta(G)$ into the graph $H$ obtained from $G$, $P_4$ and $K_{\Delta(G)}$ constructed in Proposition \ref{prop1}. We have $z(H)=max \{z(G), \Delta(G)\}+1$ and $\chi (H)=max\{\chi(G),\Delta(G)\}+1=\Delta(G)+1$. In case that $z(G)=\Delta(G)+1$ we have $z(H)=\Delta(G)+2 > \Delta(G)+1=\chi(H)$. In case that $z(G)\leq \Delta(G)$ we have clearly $z(H)\leq \chi(H)$. This completes the proof.
\end{proof}

\end{document}